\documentclass[a4paper,english,12pt]{article}
\usepackage[logonly]{trace}
\usepackage{babel}
\usepackage{amsfonts, amsmath, amssymb}
\usepackage{graphicx}
\usepackage{pstricks}
\usepackage{psfig}
\usepackage{multirow}
\usepackage{fancyhdr}
\usepackage{vmargin, fancybox}

\newcommand{\mathsym}[1]{{}}

\usepackage{theorem} \theorembodyfont{\upshape}

\newtheorem{theorem}{Theorem}[section]

\newtheorem{proposition}[theorem]{Proposition}

\newtheorem{definition}[theorem]{Definition}

\usepackage{graphicx}
\usepackage{pstricks}
\usepackage{psfig}

\title{Towards optimal DRP scheme for linear advection}

\author{Claire David *$\dag$ and Pierre Sagaut \thanks{Universit\'e Pierre et Marie Curie-Paris 6,
Institut Jean Le Rond d'Alembert, UMR CNRS 7190, Bo\^ite courrier
$n^0162$, 4 place Jussieu, 75252 Paris cedex 05, France - tel.
(+33)1.44.27.62.13; Fax (+33)1.44.27.52.59   ($\dag$ corresponding
author: { david@lmm.jussieu.fr}).}}

\begin{document}

\maketitle

\begin{abstract}

Finite difference schemes are here solved by means of a linear
matrix equation. The theoretical study of the related algebraic
system is exposed, and enables us to minimize the error due to a
finite difference approximation, while building a new DRP scheme in
the same time.
\end{abstract}

\noindent{\textbf{keywords}}\\DRP schemes, Sylvester equation

\pagestyle{myheadings} \thispagestyle{plain} \markboth{CL. DAVID AND
P. SAGAUT}{Towards optimal DRP scheme for linear advection}

\section{Introduction: Scheme classes}
\label{sec:intro}

\indent We hereafter propose a method that enables us to build a DRP
scheme while minimizing the error due to the finite difference
approximation, by means of an equivalent matrix equation.\\
\\

\noindent Consider the transport equation:
\begin{equation}
\label{transp} \frac{\partial u}{\partial t}+c\,\frac{\partial
u}{\partial x}=0 \,\,\, , \,\,\, x \,\in \,[0,L], \,\,\,t \,\in
\,[0,T]
\end{equation}

\noindent with the initial condition $u(x,t=0)=u_0(x)$.

\bigskip

\begin{proposition}

\noindent A finite difference scheme for this equation can be
written under the form:
\begin{equation} \label{scheme} {{{{{\alpha \, u}}_i}}^{n+1}}+
   {{{{{\beta \,u}}_i}}^{n}}
    +{{{{{\gamma \,u}}_i}}^{n-1}}
      +\delta \,{{{u_{i+1}}}^n}+{{{{{\varepsilon \, u}}_{i-1}}}^n}
            +{{{{{\zeta \,u}}_{i+1}}}^{n+1}}
              +{{{{{\eta \,u}}_{i-1}}}^{n-1}}+{{{{{\theta \,u}}_{i-1}}}^{n+1}}+\vartheta \,{{ u}_{i+1}}^{n-1} =0
              \end{equation}

\noindent where:
\begin{equation}
{u_l}^m=u\,(l\,h, m\,\tau)
\end{equation}
\noindent  $l\, \in \, \{i-1,\, i, \, i+1\}$, $m \, \in \, \{n-1,\,
n, \, n+1\}$, $j=0, \, ..., \, n_x$, $n=0, \, ..., \, n_t$, $h$,
$\tau$ denoting respectively the mesh size and time step ($L=n_x\,h$, $T=n_t\,\tau$).\\
The Courant-Friedrichs-Lewy number ($cfl$) is defined as $\sigma = c \,\tau / h$ .\\
\\

A numerical scheme is  specified by selecting appropriate values of
the coefficients  $\alpha$, $\beta$, $\gamma$, $\delta$,
$\varepsilon$, $\zeta$, $\eta$, $\theta$ and  $\vartheta$ in
equation (\ref{scheme}), which, for sake of usefulness, will be
written as:
\begin{equation}
\alpha=\alpha_x+\alpha_t\,\,\, , \,\,\,\beta=\beta_x+\beta_t\,\,\, ,
\,\,\,\gamma=\gamma_x+\gamma_t\,\,\, ,
\,\,\,\delta=\delta_x+\delta_t\,\,\, ,
\,\,\,\varepsilon=\varepsilon_x+\varepsilon_t\,\,\, , \,\,\,
\end{equation}
\noindent where the "$_x$" denotes a dependance upon the mesh size
$h$, while the "$_t$" denotes a dependance upon the time step
$\tau$.\\

\end{proposition}
\bigskip

\noindent The number of time steps will be denoted $n_t$, the number
of space
steps, $n_x$. In general, $n_x\gg n_t$.\\

\noindent In the following: the only dependance of the coefficients
upon the time step $\tau$ existing only in the Crank-Nicolson
scheme, we will restrain our study to the specific case:

\begin{equation}
 \alpha_t=\gamma_t=\zeta=\eta=\theta=\vartheta=0
\end{equation}

\noindent The paper is organized as follows. The building of the DRP
scheme is exposed in section \ref{DRP}. The equivalent matrix
equation, which enables us to minimize the error due to the finite
difference approximation, is presented in section \ref{Sylv}. A
numerical example is given in section \ref{Ex}.

\section{The DRP scheme}
\label{DRP}

 \noindent The first derivative $\frac{\partial
u}{\partial x}$ is approximated at the $l^{th}$ node of the spatial
mesh by:

\begin{equation}\label{approx}
 (\, \frac{\partial u}{\partial
x}\,)_l  \simeq
   {{{{{\beta_x \,u}}_{l+i}}}^{n}}
      +\delta_x \,{{{u_{l+i+1}}}^n}+{{{{{\varepsilon_x\,
      u}}_{l+i-1}}}^n}
\end{equation}
\noindent Following the method exposed by C. Tam and J. Webb in
\cite{Tam}, the coefficients $\beta_x$, $\delta_x$, and
$\varepsilon_x$ are determined requiring the Fourier Transform of
the finite difference scheme (\ref{approx}) to be a close
approximation of the partial derivative $ (\, \frac{\partial
u}{\partial x}\,
)_l$.\\
\noindent (\ref{approx}) is a special case of:

\begin{equation}\label{approx_Cont}
 (\, \frac{\partial u}{\partial
x}\,)_l  \simeq
   \beta_x \,u(x+i\,h)
      +\delta_x \,u(x+(i+1)\,h)+\varepsilon_x\,u(x+(i-1)\,h)
\end{equation}

\noindent where $x$ is a continuous variable, and can be recovered
setting $x=l\,h$.\\
\noindent Denote by $\omega$ the phase. Applying the Fourier
transform, referred to by $\,\widehat{\, }$ , to both sides of
(\ref{approx_Cont}), yields:

\begin{equation}
\label{Wavenb}
 j\, \omega \, \widehat{u}  \simeq \left \lbrace
   \beta_x \,e^{\,0}
      +\delta_x \,e^{\,j\,\omega\,h}+\varepsilon_x\,e^{\,-\,j\,\omega\,h}
     \right \rbrace \, \widehat{u}
\end{equation}
\noindent  $j$ denoting the complex square root of $-1$.\\




\noindent Comparing the two sides of (\ref{Wavenb}) enables us to
identify the wavenumber $ \overline{\lambda}$ of the finite
difference scheme (\ref{approx}) and the quantity $\frac{1}{j}\,
\left \lbrace\beta_x \,e^{\,0}
      +\delta_x
      \,e^{\,j\,\omega\,h}+\varepsilon_x\,e^{\,-\,j\,\omega\,h}\,\right
      \rbrace$, i. e.:
\noindent The wavenumber of the finite difference scheme
(\ref{approx}) is thus:

\begin{equation}
 \overline{\lambda}=-\,j\, \left \lbrace\beta_x \,e^{\,0}
      +\delta_x
      \,e^{\,j\,\omega\,h}+\varepsilon_x\,e^{\,-\,j\,\omega\,h}\,\right \rbrace
\end{equation}

\noindent To ensure that the Fourier transform of the finite
difference scheme is a good approximation of the partial derivative
$ (\, \frac{\partial u}{\partial x}\, )_l$ over the range of waves
with wavelength longer than $4\,h$, the a priori unknowns
coefficients $\beta_x$, $\delta_x$, and $\varepsilon_x$ must be
choosen so as to minimize the integrated error:

\begin{equation}\begin{array}{rcl}
 {\mathcal E} &=&\int_{-\frac{\pi}{2}}^{\frac{\pi}{2}} | \lambda \,h- \overline{\lambda}
 \,h|^2\,d(\lambda \,h)\\
 &=& \int_{-\frac{\pi}{2}}^{\frac{\pi}{2}} | \kappa+j\,h\,\left \lbrace\beta_x \,e^{\,0}
      +\delta_x
      \,e^{\,j\,\kappa}+\varepsilon_x\,e^{\,-\,j\,\kappa}\,\right \rbrace\,|^2\,d(\kappa)
      \end{array}
\end{equation}

\noindent The conditions that ${\mathcal E}$ is a minimum are:

\begin{equation}
\frac{\partial  {\mathcal E}}{\partial \beta_x}=\frac{\partial
 {\mathcal E}}{\partial \delta_x}=\frac{\partial  {\mathcal E}}{\partial \varepsilon_x}=0
\end{equation}

\noindent and provide the following system of linear algebraic
equations:

\begin{equation}\left \lbrace
\begin{array}{rcl}
2 \,\pi \,h\, \beta_x +4 \,(h\,\delta_x +h\,\varepsilon_x -1)&=&0\\
 4 \,h\,\beta_x +\pi \, (2\,
\delta_x -1)&=&0\\
4 \,h\,\beta_x +2 \,\pi \, h\,\varepsilon_x &=&0
\end{array} \right.
\end{equation}

\noindent  which enables us to determine the required values of
$\beta_x$, $\delta_x$, and $\varepsilon_x$:

\begin{equation}
\label{Opt_Values1}\left \lbrace
\begin{array}{rcl}
 \beta_x &=& \beta_x^{opt} \,=\, \frac{\pi }{h\,(\pi ^2-8)}\\
\delta_x &=&\delta_x^{opt}\,=\, \frac{1}{2}-\frac{2}{h\,(\pi ^2-8)}\\
\varepsilon_x &=&\varepsilon_x^{opt}\,=\, -\frac{2}{h\,(\pi ^2-8)}\\
\end{array} \right.
\end{equation}

\section{The Sylvester equation}
\label{Sylv}
\subsection{Matricial form of the finite differences problem}

\begin{theorem}
The problem (\ref{scheme}) can be written under the following
matricial form:
\begin{equation}
\label{Eq} {M_1}\,U +U\,M_2+{\cal{L}}(U)=M_0
\end{equation}

\noindent where $M_1$ and $M_2$ are square matrices respectively
$n_x-1$ by $n_x-1$, $n_t$ by $n_t$, given by:
\begin{equation}
\begin{array}{ccc}
{M_1}= \left (
\begin{array}{ccccc}
\beta & \delta & 0 & \ldots &  0 \\
\varepsilon& \beta & \ddots & \ddots &  \vdots \\
0 &  \ddots & \ddots & \ddots &  0\\
\vdots &  \ddots & \ddots & \beta &  \delta\\
0 &  \ldots & 0 & \varepsilon &  \beta\\
\end{array} \right ) &
 & {M_2}= \left (
\begin{array}{ccccc}
0 & \gamma& 0 & \ldots &  0 \\
\alpha & 0 & \ddots & \ddots &  \vdots \\
0 &  \ddots & \ddots & \ddots &  0\\
\vdots &  \ddots & \ddots & \ddots &  \gamma\\
0 &  \ldots & 0 & \alpha &  0\\
\end{array} \right )
\end{array}
  \end{equation}

\noindent the matrix $M_0$ being given by:
\bigskip
  \begin{equation}
  \label{M0}
\scriptsize{{M_0}= \left (
\begin{array}{ccccc}
    -\gamma \,u_1^0
      -\varepsilon \,u_{0}^1-\eta\,u_0^0-\theta\,u_{0}^{2}-\vartheta\,u_{2}^{0}
            &   -\varepsilon \,u_{0}^2 -\eta\,u_0^1-\theta\,u_{0}^{3}& \ldots & \ldots   &
            -\varepsilon \,u_{0}^{n_t}-\eta\,u_{0}^{n_t-1} \\
 -\gamma \,u_2^0-\eta\,u_{1}^{0}-\vartheta\,u_{3}^{0}
       &  0
 & \ldots  & \ldots    &  0 \\
\vdots  &  \vdots & \vdots & \vdots &  \vdots\\
 -\gamma \,u_{n_x-2}^0-\eta\,u_{n_x-2}^{0}-\vartheta\,u_{n_x-1}^{0}
       &  0
 & \ldots  & \ldots   &  0 \\
    -\gamma \,u_{n_x-1}^0-\delta \,u_{n_x}^{1}-\eta\,u_{n_x-2}^{0}-\zeta\,u_{n_x}^{2}-\vartheta\,u_{n_x}^{0}
      &   -\delta \,u_{n_x}^{2} -\zeta\,u_{n_x}^{3}-\vartheta\,u_{n_x}^{1}& \ldots
       & \ldots &  -\delta \,u_{n_x}^{n_t}-\vartheta\,u_{n_x}^{n_t-1}\\
\end{array} \right )}
  \end{equation}

\bigskip
\noindent and where ${\cal {L}}$  is a linear matricial operator
which can be written as: \begin{equation} {\cal {L}}={\cal
{L}}_1+{\cal {L}}_2+{\cal {L}}_3+{\cal {L}}_4
\end{equation}
\noindent where ${\cal {L}}_1$, ${\cal {L}}_2$, ${\cal {L}}_3$ and
${\cal {L}}_4$ are given by:
\begin{equation}
\begin{array}{ccc}
 {\cal {L}}_1(U)  = \zeta  \left (
\begin{array}{ccccc}
u_2^2 & u_2^3 & \ldots & u_2^{n_t} &  0\\
u_3^2 & u_3^3 & \ldots & \vdots &   \vdots \\
\vdots &  \vdots & \ddots & \vdots &  \vdots\\
u_{n_x-1}^2 & u_{n_x-1}^3  & \ldots &  u_{n_x-1}^{n_t}& 0  \\
0 & 0  & \ldots &  0& 0  \\
\end{array}
 \right )
  &
  &  {\cal {L}}_2(U) = \eta \left (
\begin{array}{ccccc}
0 & 0  & \ldots &  0& 0 \\
0 & u_1^1  & u_1^2 & \ldots &  u_1^{n_t-1} \\
0 & u_1^0  & u_1^1 & \ldots &  u_2^{n_t-1} \\
\vdots &  \vdots & \vdots &  \ddots & \vdots \\
0 & u_{n_x-2}^1  & u_{n_x-2}^2 & \ldots &  u_{n_x-2}^{n_t-1} \\
\end{array}
 \right )
\end{array}
  \end{equation}

\begin{equation}
\begin{array}{ccc}
{\cal {L}}_3(U) = \theta \left (
\begin{array}{ccccc}
0 & \ldots &  \ldots &    \ldots &  0 \\
u_1^2 & u_1^3  & \ldots &  u_1^{n_t}& 0 \\
u_2^2 & u_2^3  & \ldots &  u_2^{n_t}& 0 \\
\vdots &  \vdots & \vdots &  \vdots & \vdots \\
u_{n_x-2}^2 & u_{n_x-2}^3  & \ldots &  u_{n_x-2}^{n_t}& 0 \\
\end{array}
 \right )
  &
  &  {\cal {L}}_4(U) = \vartheta \left (
\begin{array}{ccccc}
0 & u_2^1 & u_2^2 & \ldots  &  u_2^{n_t-1} \\
0 & u_3^1 & u_3^2 & \ldots  &  u_3^{n_t-1} \\
\vdots &  \vdots & \ddots  & \ddots &  \vdots\\
0 & u_{n_x-1}^1& \ldots & \ldots  &  u_{n_x-1}^{n_t-1} \\
0 &  0 & \ldots &   \ldots &  0\\
\end{array}
 \right )
 \end{array}
 \end{equation}
\bigskip

\end{theorem}

\bigskip

\begin{proposition}
\noindent The second member matrix $M_0$ bears the initial
conditions, given for the specific value $n=0$, which correspond to
the initialization process when computing loops, and the boundary
conditions, given for the specific values $i=0$, $i=n_x$.
\end{proposition}

\bigskip

\noindent Denote by $u_{exact}$ the exact solution of (\ref{transp}).\\
\noindent The corresponding matrix $U_{exact}$ will be:

\begin{equation}
U_{exact}=[{{{U_{{exact}_i}}}^n}{]_{\, 1\leq i\leq {n_x-1},\, 1\leq
n\leq {n_t}\, }} \end{equation} where:

\begin{equation}
{U_{exact}}_i^n=U_{exact}(x_i,t_n)
 \end{equation}

\noindent with $x_i=i \; h$, $t_n=n \; \tau$. \bigskip

\bigskip

  \begin{definition}
\noindent We will call \textit{error matrix} the matrix defined by:
 \begin{equation}
 \label{err}
E=U-U_{exact}
   \end{equation}
  \end{definition}

\bigskip

\noindent Consider the matrix $F$ defined by:
\begin{equation} F={M_1}\,U_{exact}+U_{exact}\,M_2 + {\cal{L}}(U_{exact})-M_0\end{equation}

\bigskip

  \begin{proposition}
\noindent The \textit{error matrix} $E$ satisfies:

   \begin{equation}
   \label{eqmtr}
{M_1}\,E+E\,M_2+{\cal{L}}(E)=F
   \end{equation}

 \end{proposition}

\subsection{The matrix equation}

   \begin{theorem}
\noindent Minimizing the error due to the approximation induced by
the numerical scheme is equivalent to minimizing the norm of the
matrices $E$ satisfying (\ref{eqmtr}).
   \end{theorem}

   \bigskip

{\em Note:} \noindent Since the linear matricial operator
${\cal{L}}$ appears only in the Crank-Nicolson scheme, we will
restrain our study to the case ${\cal{L}}=0$. The generalization to
the case ${\cal{L}} \neq 0$ can be easily deduced.

   \bigskip

   \begin{proposition}

   \noindent The problem is then the determination of the minimum norm solution
of:

   \begin{equation}
   \label{SylvErr}
{M_1}\,E+E\,M_2=F
   \end{equation}

\noindent which is a specific form of the Sylvester equation:

   \begin{equation}
   \label{SylvGen}
AX+XB=C
   \end{equation}
where $A$ and $B$ are respectively $m$ by $m$ and $n$ by $n$
matrices, $C$ and $X$, $m$ by $n$ matrices.

   \end{proposition}

\bigskip

\subsection{Minimization of the error} \label{MinErr}

\noindent Calculation yields:

\footnotesize \noindent \begin{equation}\left \lbrace
\begin{array}{ccc}
{M_1}\,^T M_1&=& diag \big
 (\left ( \begin{array}{cc}
\beta^2+ \delta^2 & \beta\,(\delta+\varepsilon) \\
\beta\,(\delta+\varepsilon) &\varepsilon^2+ \beta^2 \\
\end{array} \right ),\ldots, \left ( \begin{array}{cc}
\beta^2+ \delta^2 & \beta\,(\delta+\varepsilon) \\
\beta\,(\delta+\varepsilon) &\varepsilon^2+ \beta^2 \\
\end{array} \right )
 \big )\\ {M_2}\,^T M_2&=&
 diag \big
 (\left ( \begin{array}{cc}
\gamma^2 & 0 \\
0 &\alpha^2\\
\end{array} \right ),\ldots, \left ( \begin{array}{cc}
\gamma^2 & 0 \\
0 &\alpha^2 \\
\end{array} \right )
\end{array} \right.
  \end{equation}

\normalsize \noindent The singular values of $M_1$ are the singular
values of the block matrix $\big
 (\left ( \begin{array}{cc}
\beta^2+ \delta^2 & \beta\,(\delta+\varepsilon) \\
\beta\,(\delta+\varepsilon) &\varepsilon^2+ \beta^2 \\
\end{array} \right )$, i. e. \begin{equation} \frac{1}{2} \,(2 \beta ^2+\delta
^2+\varepsilon
   ^2-(\delta +\varepsilon ) \,\sqrt{4 \beta ^2+\delta
   ^2+\varepsilon ^2-2 \delta \, \varepsilon })  \end{equation} \noindent of order $\frac {n_x-1}{2}$, and \begin{equation} \frac{1}{2} \,(2 \beta
^2+\delta ^2+\varepsilon
   ^2+(\delta +\varepsilon ) \,\sqrt{4 \beta ^2+\delta
   ^2+\varepsilon ^2-2 \delta \, \varepsilon })  \end{equation} \noindent of order $\frac {n_x-1}{2}$.\\

\noindent The singular values of $M_2$ are $\alpha^2 $, of order
$\frac {n_t}{2}$, and $\gamma^2 $, of order $\frac {n_t}{2}$.

\noindent Consider the singular value decomposition of the matrices
$M_1$ and $M_2$:

\begin{equation}
U_1^T\,M_1\,V_1=\left (
\begin{array}{cc}
\widetilde{M_1} &0 \\
0& 0\\
\end{array} \right )
\,\,\, ,  \,\,\, U_2^T\,M_1\,V_2=\left (
\begin{array}{cc}
\widetilde{M_2} &0 \\
0& 0\\
\end{array} \right )
  \end{equation}

\noindent where $U_1$, $V_1$, $U_2$, $V_2$, are orthogonal matrices.
\noindent $\widetilde{M_1}$, $\widetilde{M_2}$ are diagonal
matrices, the diagonal terms of which are respectively the nonzero
eigenvalues of the symmetric matrices $M_1\,^T M_1$, $M_2\,^T
M_2$.\\

\noindent Multiplying respectively \ref{SylvErr} on the left side by
$^T U_1$, on the right side by $V_2$, yields:

\begin{equation}
U_1^T\,M_1\,E\,V_2+U_1^T\,E\,M_2\,V_2=U_1^T\,F\,V_2
  \end{equation}

\noindent which can also be taken as:
\begin{equation}
^T U_1\,M_1\,V_1 \,^T V_1\,E\,V_2+^T U_1\,E\,^T U_2\,^T
U_2\,M_2\,V_2=U_1^T\,F\,V_2
  \end{equation}

\noindent Set:

\begin{equation}
^T V_1\,E\,V_2= \left (
\begin{array}{cc}
\widetilde{E_{11}} & \widetilde{E_{12}} \\
\widetilde{E_{21}} & \widetilde{E_{22}}\\
\end{array} \right )
\,,\,^T U_1\,E\,^T U_2= \left (\begin{array}{cc}
\widetilde{\widetilde{E_{11}}} & \widetilde{\widetilde{E_{12}}} \\
\widetilde{\widetilde{E_{21}}} & \widetilde{\widetilde{E_{22}}}\\
\end{array} \right )
  \end{equation}

\begin{equation}
\label{Ftilde} ^T U_1\,F\,V_2= \left (
\begin{array}{cc}
\widetilde{F_{11}} & \widetilde{F_{12}} \\
\widetilde{F_{21}} & \widetilde{F_{22}}\\
\end{array} \right )
  \end{equation}
\noindent We have thus:
\begin{equation}\left (\begin{array}{cc}
\widetilde{M_{1}}\,\widetilde{E_{11}} & \widetilde{M_{1}}\,\widetilde{E_{12}} \\
0 & 0\\
\end{array} \right )+\left (\begin{array}{cc}
\widetilde{\widetilde{E_{11}}}\,\widetilde{M_{2}} & 0 \\
\widetilde{\widetilde{E_{21}}}\,\widetilde{M_{2}} & 0\\
\end{array} \right )=\left (\begin{array}{cc}
\widetilde{F_{11}} & \widetilde{F_{12}} \\
\widetilde{F_{21}} & \widetilde{F_{22}}\\
\end{array} \right )
  \end{equation}

\noindent It yields:
\begin{equation}
\left \lbrace
\begin{array}{ccc}
\widetilde{M_{1}}\,\widetilde{E_{11}} +
\widetilde{\widetilde{E_{11}}}\,\widetilde{M_{2}}&=&\widetilde{F_{11}}\\
\widetilde{M_{1}} \,\widetilde{E_{12}}&=&\widetilde{F_{12}}\\
\widetilde{\widetilde{E_{21}}}\,\widetilde{M_{2}}&=&\widetilde{F_{21}}\\
\end{array} \right.
  \end{equation}

\noindent One easily deduces:
\begin{equation}
\left \lbrace
\begin{array}{ccc}
\widetilde{E_{12}}&=&{\widetilde{M}_{1}}^{-1}\,\widetilde{F_{12}}\\
{\widetilde{\widetilde{E}_{21}}}&=&\widetilde{F_{21}}\,{\widetilde{M_{2}}}^{-1}\\
\end{array} \right.
  \end{equation}

\noindent The problem is then the determination of the
$\widetilde{E_{11}}$ and $\widetilde{\widetilde{E_{11}}}$
satisfying:

\begin{equation}
\label{Pb} \widetilde{M_{1}}\,\widetilde{E_{11}} +
\widetilde{\widetilde{E_{11}}}\,\widetilde{M_{2}}=\widetilde{F_{11}}
  \end{equation}

\noindent Denote respectively by $\widetilde{e_{ij}}$,
$\widetilde{\widetilde{e_{ij}}}$ the components of the matrices
$\widetilde{E}$, $\widetilde{\widetilde{E}}$.\\
\noindent The problem \ref{Pb} uncouples into the independent
problems:\\ \noindent minimize
\begin{equation}
\sum_{i,j} {\widetilde{e_{ij}}}^2+{\widetilde{\widetilde{e_{ij}}}}^2
\end{equation}

 \noindent under the constraint \begin{equation}
\widetilde{M_{1}}_{ii}\,
{\widetilde{e_{ij}}}+\widetilde{M_{2_{ii}}}\,{\widetilde{\widetilde{e_{ij}}}}=\widetilde{F_{11}}_{ij}
  \end{equation}
\noindent This latter problem has the solution:

\begin{equation}\left \lbrace
\begin{array}{ccc}
\widetilde{e_{ij}}
&=&\frac{\widetilde{{M_{1}}_{ii}}\,\widetilde{{F_{11}}_{ij}}}
{{\widetilde{{M_{1}}_{ii}}}^2+{\widetilde{{M_{2}}_{jj}}^2}}\\
\widetilde{\widetilde{e_{ij}}} &=&\frac{\widetilde{{M_{2}}
_{jj}}\,\widetilde{{F_{11}}_{ij}}}{{\widetilde{{M_{1}}_{ii}}}^2+{\widetilde{{M_{2}}_{jj}}^2}}\\
\end{array} \right.
  \end{equation}

\noindent The minimum norm solution of \ref{SylvErr} will then be
obtained when the norm of the matrix $\widetilde{{F_{11}}}$ is
minimum.\\
\noindent In the following, the euclidean norm will be considered.

\noindent Due to (\ref{Ftilde}):
 \begin{equation}
  \|\widetilde{{F_{11}}}  \| \leq  \|\widetilde{{F}}  \| \leq \|U_1
  \| \,\|F \|\,\|V_2  \| \leq \|U_1
  \| \,\|V_2  \| \, \|M_1\,U_{exact}+U_{exact}\,M_2-M_0
  \|
   \end{equation}

\noindent $U_1$ and $V_2$ being orthogonal matrices, respectively
$n_x-1$ by $n_x-1$, $n_t$ by $n_t$, we have:

 \begin{equation}
  \|U_1
  \|^2 =n_x-1\,\,\,, \,\,\,\|V_2
  \|^2 =n_t
   \end{equation}

\noindent Also:

 \begin{equation}
  \|M_1
  \|^2 =\frac{n_x-1}{2}\,\big ( 2\,\beta^2+\delta^2+\varepsilon^2 \big )\,\,\,, \,\,\,
  \|M_2
  \|^2 =\frac{n_t}{2}\,\big ( \alpha^2+\gamma^2 \big )
   \end{equation}

\noindent The norm of $M_0$ is obtained thanks to relation
(\ref{M0}).

\noindent This results in:
 \begin{equation}
 \label{Min}
  \|\widetilde{{F_{11}}} \| \leq \sqrt {n_t\,(n_x-1)} \, \left \lbrace  \| U_{exact} \| \,\big (
  \sqrt{\frac{n_x-1}{2}}\,\sqrt{ 2\,\beta^2+\delta^2+\varepsilon^2 }+
   \sqrt{\frac{n_t}{2}}\,\sqrt{\alpha^2+\gamma^2 } \,\big )+
  \|M_0  \|\right \rbrace
   \end{equation}

\noindent $ \|\widetilde{{F_{11}}} \| $ can be minimized through the
minimization of the second factor of the right-side member of
(\ref{Min}), which is function of the scheme parameters.\\

\noindent $\| U_{exact} \| $ is a constant. The quantities $
\sqrt{\frac{n_x-1}{2}}\,\sqrt{ 2\,\beta^2+\delta^2+\varepsilon^2
  }$, $\sqrt{\alpha^2+\gamma^2 } $ and $\|M_0  \|$ being
  strictly positive, minimizing the second factor of the right-side member
  of (\ref{Min}) can be obtained through the minimization of the
  following functions:

 \begin{equation}
 \label{func_Sylv}
 \left \lbrace
\begin{array}{rcl}
f_1(\beta,\delta,\varepsilon)&=&\sqrt{ 2\,\beta^2+\delta^2+\varepsilon^2}\\
 f_2(\alpha,\gamma)&=&\sqrt{\alpha^2+\gamma^2 }\\
 f_3(\alpha,\beta,\gamma,\delta,\varepsilon)&=&\|M_0  \|\\
\end{array} \right.
\end{equation}

\noindent i.e.:

 \begin{equation}
 \left \lbrace
\begin{array}{rcl}
f_1(\beta,\delta,\varepsilon)&=&\sqrt{ 2\,(\beta_x+\beta_t)^2+(\delta_x+\delta_t)^2+(\varepsilon_x+\varepsilon_t)^2}\\
 f_2(\alpha,\gamma)&=&\sqrt{\alpha_x^2+\gamma_x^2 }\\
 f_3(\alpha,\beta,\gamma,\delta,\varepsilon)&=&\|M_0  \|\\
\end{array} \right.
\end{equation}

\noindent Setting:

\begin{equation}
\left \lbrace
\begin{array}{rcl}
 \beta_x &=& \beta_x^{opt} \\
\delta_x &=&\delta_x^{opt}\\
\varepsilon_x &=&\varepsilon_x^{opt}\\
\end{array} \right.
\end{equation}

\noindent one obtains the DRP scheme with the minimal error through
the minimization of:

 \begin{equation}
 \left \lbrace
\begin{array}{rcl}
g_1(\beta_t,\delta_t,\varepsilon_t)&=&\sqrt{
\,(\beta_x^{opt}+\beta_t)^2+(\delta_x^{opt}+\delta_t)^2
+(\varepsilon_x^{opt}+\varepsilon_t)^2}\\
 g_2(\alpha_x,\gamma_x)&=&\sqrt{\alpha_x^2+\gamma_x^2 }\\
g_3(\alpha,\beta_t,\gamma,\delta_t,\varepsilon_t)&=&\|M_0  \|\\
\end{array} \right.
\end{equation}

\section{Numerical results}

\begin{figure}[h!]
\center{\includegraphics[height=5cm]{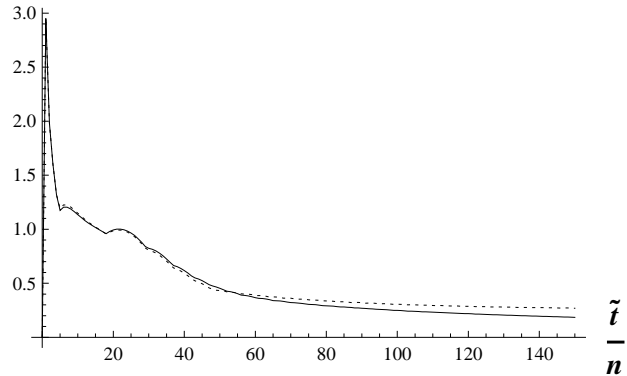}}\\
\caption{\small{$L^\infty$ norm of the error for the optimized
scheme (in black) and the non optimized one (dashed curve)}.}
\label{Erreur1}
\end{figure}

\noindent We denote by $\tilde{t}$ the non-dimensional time
parameter. Figure \ref{Erreur1} displays the $L^\infty$ norm of
the error for an optimized scheme (in black), where $\beta_x$,
$\delta_x$, $\varepsilon_x$ are given by (\ref{Opt_Values1}), and
a non-optimized one: numerical results perfectly fit the
theoretical ones.

\begin{figure}[h!]
\center{\includegraphics[height=5cm]{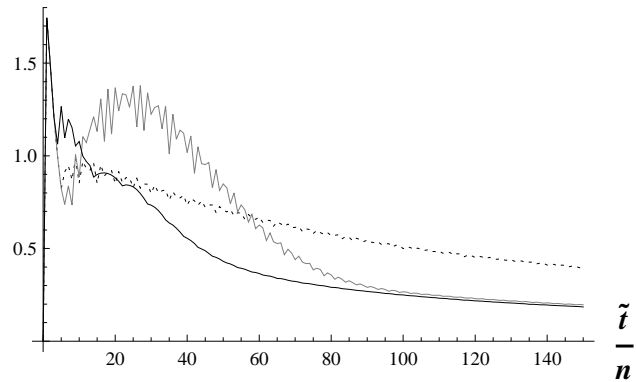}}\\
\caption{\small{$L^\infty$ norm of the error for the optimized
scheme (in black), the DRP scheme (in gray), and the FCTS scheme
(dashed plot)}.} \label{Erreur2}
\end{figure}

\noindent Figure \ref{Erreur2} displays the $L^\infty$ norm of the
error for the above optimized scheme  (in black), a seven-point
stencil DRP scheme (in gray), and the FCTS scheme (dashed plot).
As time increases, the optimized scheme yields, as expected,
better results than the FCTS one. Also, for $15 \leq \frac
{\tilde{t}} {n} \leq 100$, results appear to be better than those
of the classical DRP scheme. For large values of the time
parameter, both latter schemes yield the same results.

\begin{figure}[h!]
\center{\includegraphics[height=5cm]{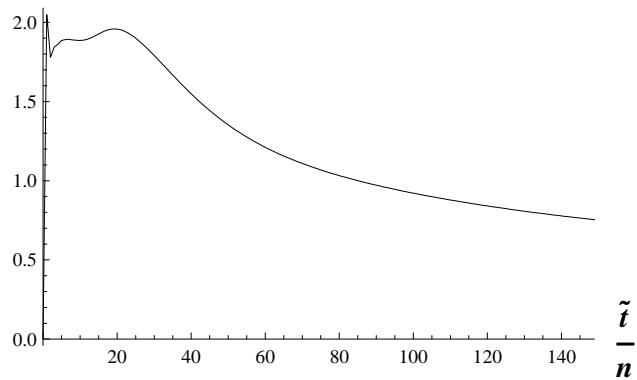}}\\
\caption{\small{$L^2$ norm of the error for the optimized case (in
black) and the optimized case (dashed curve)}.} \label{ErreurL2}
\end{figure}

\noindent Figure \ref{ErreurL2} displays the $L^2$ norm of the
error for the above optimized scheme (in black), the seventh-order
DRP scheme (in gray), and the FCTS scheme (dashed plot). As
expected, results coincide.


\section{Conclusion}

The above results open new ways for the building of DRP schemes. It
seems that the research on this problem has not been performed
before as far as our knowledge goes. In the near future, we are
going to extend the techniques described herein to nonlinear
schemes, in conjunction with other innovative methods as the Lie
group theory.

\addcontentsline{toc}{section}{\numberline{}References}

\end{document}